\documentclass[12pt]{article}

\usepackage{amsfonts}

\newtheorem{thm}{Theorem}[section]
\newtheorem{lemma}[thm]{Lemma}
\newtheorem{cor}[thm]{Corollary}

\newenvironment{remark}{\par\medskip\noindent{\bf Remark.\ }}{\par\smallskip}
\newcommand{\proof
}{\par\medskip\noindent {\bf Proof.\ \ }}

\newcommand{\be}{\begin{equation}}
\newcommand{\ee}{\end{equation}}
\newcommand{\openbox}{\leavevmode
  \hbox to8pt{\hfil\vrule\vbox to6pt{\hrule width6pt\vfil\hrule}\vrule}}

\newcommand{\qed}{\hbox to5pt{ } \hfill \openbox\bigskip\medskip}

\newcommand{\Z}{\mathbb Z}
\newcommand{\Q}{\mathbb Q}

\title{Betti numbers of Stanley--Reisner rings with pure resolutions}
\author{G\'abor Heged\"{u}s
\\{\normalsize Johann Radon Institute for Computational and Applied Mathematics}
}

\begin{document}

\footnotetext{
{\bf Keywords.}  Betti number, Hilbert function, Stanley-Reisner ring

{\bf 2000 Mathematics Subject Classification.}  05E40, 13D02, 13D40 }

\maketitle

\begin{abstract}
Let $\Delta$ be simplicial complex and let $k[\Delta]$ denote the 
Stanley--Reisner ring corresponding to $\Delta$. 
Suppose that $k[\Delta]$ has a pure free resolution. 
Then we describe the Betti numbers and the Hilbert--Samuel multiplicity 
of $k[\Delta]$ in terms of the $h$--vector 
of $\Delta$. As an application, we derive a linear equation system 
for the components of 
the $h$--vector
of the clique complex of an arbitrary chordal graph.

\end{abstract}
\medskip


\section{Introduction}
\noindent

Let $k$ denote an arbitrary field.
Let $R$ be the graded ring $k[x_1,\ldots ,x_n]$. 
The vector space $R_s=k[x_1,\ldots ,x_n]_s$ consists of 
the homogeneous polynomials of total degree $s$, together with $0$.

In \cite{Fr} R. Fr\"oberg characterized the graphs $G$ such that $G$ 
has a linear free resolution. He proved:

\begin{thm} \label{Froberg_theorem}
Let $G$ be  a simple graph on $n$ vertices. Then $R/I(G)$ has linear 
free resolution precisely when $\overline{G}$, the complementary graph of $G$ 
is chordal. 
\end{thm}

In \cite{E2} E. Emtander generalized Theorem \ref{Froberg_theorem} for 
generalized chordal hypergraphs.

In this article we prove explicit formulas for the Betti numbers 
of the Stanley--Reisner ring $k[\Delta]$
such that $k[\Delta]$ has a pure free resolution in terms of the $h$--vector
of $\Delta$.  

In Section 2 we collected  some basic results about simplicial complices,
 free resolutions, Hilbert fuctions and Hilbert series.
 We present our main results in Section 3.


\section{Preliminaries}
\subsection{Free resolutions}

Recall that for every finitely generated graded module $M$ over $R$ 
we can associate to $M$ 
a {\em minimal graded free resolution} 
$$ 
0\longrightarrow \bigoplus_{i=1}^{\beta_p} R(-d_{p,i}) \longrightarrow 
\bigoplus_{i=1}^{\beta_{p-1}} R(-d_{p-1,i})\longrightarrow \ldots \longrightarrow 
\bigoplus_{i=1}^{\beta_{0}} R(-d_{0,i}) \longrightarrow M \longrightarrow
0, 
$$
where $p\leq n$ and $R(-j)$ is the free $R$-module obtained 
by shifting the degrees of $R$ by $j$. 

Here the natural number ${\beta}_{k}$ is the $k$'th {\em total
Betti number} of $M$ and $p$ is the projective dimension of $M$. 

The module $M$ has a {\em pure resolution} if there are 
constants $d_0<\ldots < d_p$
such that 
$$
d_{0,i}=d_0,\ldots ,d_{p,i}=d_p 
$$
for all $i$. If in addition
$$
d_i=d_0+i,
$$
for all $1\leq i\leq p$, then 
we call the minimal free resolution  to be {\em $d_0$--linear}. 

In \cite{R} Theorem 2.7 the following bound for the Betti numbers was proved.
\begin{thm} \label{Betti_bound}
Let $M$ be an $R$--module having a pure resolution
of type $(d_0,\ldots, d_p)$ and Betti numbers $\beta_0, \ldots , \beta_p$, where 
$p$ is the projective dimension of $M$.
Then 
\begin{equation} \label{bound_Betti}
\beta_i\geq {p\choose i}
\end{equation}
for each $0\leq i\leq p$.
\end{thm}

\subsection{Hilbert--Serre Theorem}

Let $M=\bigoplus_{i\geq 0} M_i$ be a finitely generated 
nonnegatively graded module over the polynomial ring $R$.
We call the formal power series 
$$
H_M(z):= \sum_{i=0}^{\infty} h_M(i)z^i
$$
the {\em Hilbert--series} of the module $M$.

The Theorem of Hilbert--Serre states that there exists 
a (unique) polynomial $P_M(z)\in \Q[z]$, the so-called {\em Hilbert polynomial}
of $M$, such that $h_M(i)=P_M(i)$ for each $i>>0$. Moreover, $P_M$ has degree
$\mbox{dim }M-1$ and $(\mbox{dim }M-1)!$ times the leading coefficient 
of $P_M$ is the {\em Hilbert--Samuel multiplicity} of $M$, denoted here 
by $e(M)$. 

Hence there exist integers $m_0,\ldots, m_{d-1}$ 
such that $h_M(z)=m_0\cdot{z\choose d-1}+m_1\cdot{z\choose d-2}+\ldots + m_{d-1}$,
where ${z\choose r}=\frac{1}{r!}z(z-1)\ldots (z-r+1)$ and $d:=\mbox{dim} M$. 
Clearly $m_0=e(M)$.

We can summarize the Hilbert-Serre theorem as follows: 
\begin{thm} (Hilbert--Serre) \label{Hilbert_Serre}
Let $M$ be a finitely generated nonnegatively graded 
$R$--module of dimension $d$, then the following stetements hold:\\
(a) There exists a (unique) polynomial $P(z)\in \Z[z]$ such that 
the Hilbert--series $H_M(z)$ of $M$ may be written as 
$$
H_M(z)=\frac{P(z)}{(1-z)^d}
$$
(b) $d$ is the least integer for which $(1-z)^dH_M(z)$ is a polynomial.
\end{thm}


\subsection{Simplicial complices and Stanley--Reisner rings}

We say that $\Delta\subseteq 2^{[n]}$ is a {\em simplicial complex}
 on the vertex set $[n]=\{1,2,\ldots ,n\}$, if 
 $\Delta$ is a set of subsets of $[n]$ such that 
$\Delta$ is a down--set, that is, $G\in\Delta$ and $F\subseteq G$ 
implies that $F\in \Delta$, and $\{i\}\in \Delta$ for all $i$.

The elements of $\Delta$ are called {\em faces} 
and the {\em dimension} of a face is one less than its cardinality. An $r$-face is an abbreviation for an $r$-dimensional face.
The dimension of $\Delta$ is the dimension of a maximal face.
We use the notation $\mbox{dim}(\Delta)$ for the dimension 
of $\Delta$.

If $\mbox{dim}(\Delta)=d-1$, then 
the $(d+1)$--tuple $(f_{-1}(\Delta),\ldots ,f_{d-1}(\Delta))$ 
is called the {\em $f$-vector} of $\Delta$, where
 $f_i(\Delta)$ denotes the number of $i$--dimensional faces 
of $\Delta$.

Let $\Delta$ be an arbitrary simplicial complex
on $[n]$. The {\em Stanley--Reisner ring} 
$k[\Delta]:=R/I(\Delta)$ of $\Delta$
is the quotient of the ring $R$ by the
 {\em Stanley--Reisner ideal}
$$
I(\Delta):=\langle x^F:~ F \notin \Delta \rangle, 
$$ 
generated by the non--faces of $\Delta$.

The following Theorem was proved in \cite{BH} Theorem 5.1.7.
\begin{thm} \label{H-S}
Let $\Delta$ be a $d-1$--dimensional simplicial complex with $f$-vector 
$f(\Delta):=(f_{-1},\ldots ,f_{d-1})$. 
Then the Hilbert--series 
of the Stanley--Reisner ring $k[\Delta]$ is
$$
H_{k[\Delta]}(z)=\sum_{i=-1}^{d-1} \frac{f_it^{i+1}}{(1-t)^{i+1}}.
$$
\end{thm}
\qed

Recall from Theorem \ref{Hilbert_Serre}
that a homogeneous $k$-algebra $M$ of dimension $d$ has 
a Hilbert series of the forn 
$$
H_M(z)=\frac{P(z)}{(1-z)^d}
$$
where $P(z)\in \Z[z]$. 
Let $\Delta$ be a $(d-1)$--dimensional simplicial complex and write 
\begin{equation} \label{Hilbert10}
H_{k[\Delta]}(z)=\frac{\sum_{i=0}^d h_iz^i}{(1-z)^d}.
\end{equation}

\begin{lemma}
The $f$-vector and the $h$-vector of a $(d-1)$--dimensional simplicial complex
$\Delta$ are related by
$$
\sum_{i} h_it^i=\sum_{i=0}^d f_{i-1}t^i(1-t)^{d-i}.
$$
In particular, the $h$-vector has length at most $d$, and 
$$
h_j=\sum_{i=0}^j (-1)^{j-i} {d-i \choose j-i}f_{i-1}
$$
for each $j=0,\ldots ,d$.
\end{lemma}
\qed


\section{Our main result}

In the following Theorem we describe 
the Betti numbers of $k[\Delta]$ in terms of the $h$--vector 
of $\Delta$.
\begin{thm} \label{main_Betti}
Let $\Delta$ be a $(d-1)$--dimensional simplicial complex. Suppose 
that the Stanley--Reisner ring $k[\Delta]$ has a pure free resolution
\begin{equation}
{{\cal F}}_{\Delta}: 0\longrightarrow R(-d_p)^{\beta_p} \longrightarrow \ldots \longrightarrow 
\end{equation}

\begin{equation} \label{free6}
\longrightarrow R(-d_1)^{\beta_1} \longrightarrow R(-d_0)^{\beta_0} \longrightarrow 
R \longrightarrow k[\Delta] \longrightarrow 0.         
\end{equation}
Here $p$ is the projective dimension of the Stanley--Reisner ring $k[\Delta]$.

If $h(\Delta):=(h_{0}(\Delta),\ldots ,h_{d}(\Delta))$ is the 
$h$-vector of the complex $\Delta$, then
$$
\beta_i=\sum_{\ell=0}^{d_i} (-1)^{\ell+i+1}{n-d \choose \ell}h_{d_i-\ell}
$$
for each $0\leq i\leq p$.
\end{thm}
\begin{remark}
Clearly $h_i=0$ for each $i>d$.
\end{remark}

\begin{remark}
J. Herzog and M. K\"uhl proved similar formulas for the Betti number in \cite{HK} 
 Theorem 1. Here we did not assume that the Stanley--Reisner ring $k[\Delta]$
 with pure resolution is Cohen--Macaulay.   
\end{remark}
\proof
Let $M:=k[\Delta]$ denote the Stanley--Reisner ring of $\Delta$. 
Then we infer from Theorem \ref{H-S} that
\begin{equation} \label{Hilbert8}
H_M(z)=\frac{\sum_{i=0}^d h_iz^i}{(1-z)^d}.
\end{equation}
Since the Hilbert--series is additive on short exact sequences,
and since 
$$
H_R(z)=\frac{1}{(1-z)^n},
$$
and consequently
$$
H_{R(-s)}(z)=\frac{z^s}{(1-z)^n},
$$
the pure resolution
\begin{equation}
{{\cal F}}_{\Delta}: 0\longrightarrow R(-d_p)^{\beta_p} \longrightarrow \ldots \longrightarrow 
\end{equation}

\begin{equation} \label{free5}
\longrightarrow R(-d_1)^{\beta_1} \longrightarrow 
R(-d_0)^{\beta_0} \longrightarrow R \longrightarrow M \longrightarrow 0.         
\end{equation}
yields to
\begin{equation} \label{Hilbert7}
H_M(z)=\frac{1}{(1-z)^n}+\sum_{i=0}^p (-1)^{i+1} \beta_i \frac{z^{d_i}}{(1-z)^n}, 
\end{equation}
where $p=pdim(M)$.

Write $d:=\mbox{dim} M$, and let $m:=\mbox{codim}(M)=n-d$. 
It follows from the Auslander--Buchbaum formula that $m\leq p$. 

Comparing the two expressions (\ref{Hilbert7}) and (\ref{Hilbert8})
for $H_M$, we find
\begin{equation} \label{Hilbert6}
(1-z)^m \left(\sum_{i=0}^d h_iz^i \right)=\sum_{i=0}^p (-1)^{i+1} \beta_i z^{d_i}+1
\end{equation}

Using the binomial Theorem we get that
\begin{equation} \label{Hilbert9}
\left( \sum_{j=0}^{n-d} (-1)^j{n-d\choose j}z^j \right)\left(\sum_{i=0}^d h_iz^i \right)
=\sum_{i=0}^p (-1)^{i+1} \beta_i z^{d_i}
\end{equation}

Comparing the coefficients on the two sides of (\ref{Hilbert9}), we get the result.
\qed

\begin{cor}
Let $\Delta$ be a $(d-1)$--dimensional simplicial complex. Then 
$$
e(k[\Delta])=f_{d-1}.
$$
\end{cor}
\proof 
It follows from \cite{BH} Proposition 4.1.9 and (\ref{Hilbert10}) that
$$
e(k[\Delta])=\left(\sum_{i=0}^d h_iz^i\right)\mid_{z=1}=\sum_{i=0}^d h_i=f_{d-1}.
$$

\begin{cor} \label{main_cor}
Let $\Delta$ be a $(d-1)$--dimensional simplicial complex. Suppose 
that the Stanley--Reisner ring $k[\Delta]$ has an $t$--linear free resolution
\begin{equation}
{{\cal F}}_{\Delta}: 0\longrightarrow R(-t-p)^{\beta_p} \longrightarrow \ldots \longrightarrow 
\end{equation}

\begin{equation} \label{free7}
\longrightarrow R(-t-1)^{\beta_1} \longrightarrow R(-t)^{\beta_0} 
\longrightarrow R \longrightarrow k[\Delta] \longrightarrow 0.         
\end{equation}
Here $p$ is the projective dimension of the Stanley--Reisner ring $k[\Delta]$.

If $h(\Delta):=(h_{0}(\Delta),\ldots ,h_{d}(\Delta))$ is the 
$h$-vector of the complex $\Delta$, then
$$
\beta_i=\sum_{\ell=0}^{t+i} (-1)^{\ell+i+1}h_{t+i-\ell}{n-d \choose \ell}
$$
for each $0\leq i\leq p$.
\end{cor}

\begin{cor} \label{main_cor3}
Let $\Delta$ be a $(d-1)$--dimensional simplicial complex. Suppose 
that the Stanley--Reisner ring $k[\Delta]$ has an $t$--linear free resolution
\begin{equation}
{{\cal F}}_{\Delta}: 0\longrightarrow R(-t-p)^{\beta_p} \longrightarrow \ldots \longrightarrow 
\end{equation}

\begin{equation} \label{free12}
\longrightarrow R(-t-1)^{\beta_1} \longrightarrow R(-t)^{\beta_0} 
\longrightarrow R \longrightarrow k[\Delta] \longrightarrow 0.         
\end{equation}
Here $p$ is the projective dimension of the Stanley--Reisner ring $k[\Delta]$.

If $h(\Delta):=(h_{0}(\Delta),\ldots ,h_{d}(\Delta))$ is the 
$h$-vector of the complex $\Delta$, then
$$
\sum_{\ell=0}^{j} (-1)^{\ell}h_{j-\ell}{n-d \choose \ell}=0.
$$
for each $j>p+t$.
\end{cor}
\proof
Let
$$
P(z):=1+\sum_{i=0}^p (-1)^{i+1} \beta_i z^{t+i}
$$
Clearly $\mbox{deg}(P)\leq p+t$. Comparing the
coefficients of both side of (\ref{Hilbert9}), we get the result. \qed

\begin{cor} \label{main_gr} 
Let $G$ be an arbitrary chordal graph. Let $\Delta:=\Delta(G)$ be the clique
 complex of
$G$ and $d:=\mbox{dim}(\Delta)+1$. 
Let $h(\Delta):=(h_{0}(\Delta),\ldots ,h_{d}(\Delta))$ denote the 
$h$-vector of the complex $\Delta$. Let $p$ be the projective 
dimension of the Stanley--Reisner ring $k[\Delta]$. Then
$$
\sum_{\ell=0}^{j} (-1)^{\ell}h_{j-\ell}{n-d \choose \ell}=0
$$
for each $j>p+2$.
\end{cor}
\proof 
This follows easily from Theorem \ref{Froberg_theorem} and 
Corollary \ref{main_cor3}.
\qed

\begin{cor} \label{main_cor2}
Let $\Delta$ be a $(d-1)$--dimensional simplicial complex. Suppose 
that the Stanley--Reisner ring $k[\Delta]$ has a pure free resolution
\begin{equation}
{{\cal F}}_{\Delta}: 0\longrightarrow R(-d_p)^{\beta_p} \longrightarrow \ldots \longrightarrow 
\end{equation}

\begin{equation} \label{free9}
\longrightarrow R(-d_1)^{\beta_1} \longrightarrow R(-d_0)^{\beta_0} 
\longrightarrow R \longrightarrow k[\Delta] \longrightarrow 0.         
\end{equation}
Here $p$ is the projective dimension of the Stanley--Reisner ring $k[\Delta]$.

Then
\begin{equation}
\sum_{\ell=0}^{d_i} (-1)^{\ell+i+1}{n-d \choose \ell}h_{d_i-\ell}\geq {p\choose i}
\end{equation}
for each $0\leq i\leq p$.
\end{cor}
\proof 
This follows easily from Theorem \ref{Betti_bound} and Theorem \ref{main_Betti}.
\qed

{\bf Acknowledgements.}  I am indebted to Josef Schicho, Russ Woodroofe 
 and Lajos R\'onyai 
for their useful remarks.

\end{document}